\documentclass[12pt]{amsart}
\usepackage{amsmath,amssymb,amsbsy,amsfonts,amsthm,latexsym,
            amsopn,amstext,amsxtra,euscript,amscd,
            mathrsfs,
            array}

\usepackage[colorlinks,linkcolor=blue,anchorcolor=blue,citecolor=blue]{hyperref}
\usepackage{color}
\begin{document}

\newtheorem{theorem}{Theorem}
\newtheorem{lemma}[theorem]{Lemma}
\newtheorem{claim}[theorem]{Claim}
\newtheorem{cor}[theorem]{Corollary}
\newtheorem{prop}[theorem]{Proposition}
\newtheorem{definition}[theorem]{Definition}
\newtheorem{question}[theorem]{Question}
\newcommand{\hh}{{{\mathrm h}}}

\numberwithin{equation}{section}
\numberwithin{theorem}{section}

\def\sssum{\mathop{\sum\!\sum\!\sum}}
\def\ssum{\mathop{\sum\ldots \sum}}

\def \balpha{\boldsymbol\alpha}
\def \bbeta{\boldsymbol\beta}
\def \bgamma{{\boldsymbol\gamma}}
\def \bomega{\boldsymbol\omega}

\def\sssum{\mathop{\sum\!\sum\!\sum}}
\def\ssum{\mathop{\sum\ldots \sum}}
\def\dsum{\mathop{\sum\  \sum}}
\def\iint{\mathop{\int\ldots \int}}

\def\squareforqed{\hbox{\rlap{$\sqcap$}$\sqcup$}}
\def\qed{\ifmmode\squareforqed\else{\unskip\nobreak\hfil
\penalty50\hskip1em\null\nobreak\hfil\squareforqed
\parfillskip=0pt\finalhyphendemerits=0\endgraf}\fi}

\newfont{\teneufm}{eufm10}
\newfont{\seveneufm}{eufm7}
\newfont{\fiveeufm}{eufm5}
%
%
\newfam\eufmfam
     \textfont\eufmfam=\teneufm
\scriptfont\eufmfam=\seveneufm
     \scriptscriptfont\eufmfam=\fiveeufm
%
%
\def\frak#1{{\fam\eufmfam\relax#1}}

\def\fK{\mathfrak K}
\def\fT{\mathfrak{T}}

\def\fA{{\mathfrak A}}
\def\fB{{\mathfrak B}}
\def\fC{{\mathfrak C}}
\def\fD{{\mathfrak D}}

\newcommand{\sX}{\ensuremath{\mathscr{X}}}

\def\eqref#1{(\ref{#1})}

\def\vec#1{\mathbf{#1}}
\def\dist{\mathrm{dist}}
\def\vol#1{\mathrm{vol}\,{#1}}

\def\squareforqed{\hbox{\rlap{$\sqcap$}$\sqcup$}}
\def\qed{\ifmmode\squareforqed\else{\unskip\nobreak\hfil
\penalty50\hskip1em\null\nobreak\hfil\squareforqed
\parfillskip=0pt\finalhyphendemerits=0\endgraf}\fi}

\def\sA{\mathscr A}
\def\sB{\mathscr B}
\def\sC{\mathscr C}
\def\sD{\Delta}
\def\sE{\mathscr E}
\def\sF{\mathscr F}
\def\sG{\mathscr G}
\def\sH{\mathscr H}
\def\sI{\mathscr I}
\def\sJ{\mathscr J}
\def\sK{\mathscr K}
\def\sL{\mathscr L}
\def\sM{\mathscr M}
\def\sN{\mathscr N}
\def\sO{\mathscr O}
\def\sP{\mathscr P}
\def\sQ{\mathscr Q}
\def\sR{\mathscr R}
\def\sS{\mathscr S}
\def\sU{\mathscr U}
\def\sT{\mathscr T}
\def\sV{\mathscr V}
\def\sW{\mathscr W}
\def\sX{\mathscr X}
\def\sY{\mathscr Y}
\def\sZ{\mathscr Z}

\def\cA{{\mathcal A}}
\def\cB{{\mathcal B}}
\def\cC{{\mathcal C}}
\def\cD{{\mathcal D}}
\def\cE{{\mathcal E}}
\def\cF{{\mathcal F}}
\def\cG{{\mathcal G}}
\def\cH{{\mathcal H}}
\def\cI{{\mathcal I}}
\def\cJ{{\mathcal J}}
\def\cK{{\mathcal K}}
\def\cL{{\mathcal L}}
\def\cM{{\mathcal M}}
\def\cN{{\mathcal N}}
\def\cO{{\mathcal O}}
\def\cP{{\mathcal P}}
\def\cQ{{\mathcal Q}}
\def\cR{{\mathcal R}}
\def\cS{{\mathcal S}}
\def\cT{{\mathcal T}}
\def\cU{{\mathcal U}}
\def\cV{{\mathcal V}}
\def\cW{{\mathcal W}}
\def\cX{{\mathcal X}}
\def\cY{{\mathcal Y}}
\def\cZ{{\mathcal Z}}
\newcommand{\rmod}[1]{\: \mbox{mod} \: #1}

\def\vr{\mathbf r}

\def\e{{\mathbf{\,e}}}
\def\ep{{\mathbf{\,e}}_p}
\def\em{{\mathbf{\,e}}_m}
\def\en{{\mathbf{\,e}}_n}

\def\Tr{{\mathrm{Tr}}}
\def\Nm{{\mathrm{Nm}}}

 \def\SS{{\mathbf{S}}}

\def\lcm{{\mathrm{lcm}}}

\def\({\left(}
\def\){\right)}
\def\fl#1{\left\lfloor#1\right\rfloor}
\def\rf#1{\left\lceil#1\right\rceil}

\def\mand{\qquad \mbox{and} \qquad}

\newcommand{\commG}[1]{\marginpar{%
\begin{color}{red}
\vskip-\baselineskip 
\raggedright\footnotesize
\itshape\hrule \smallskip G: #1\par\smallskip\hrule\end{color}}}

\newcommand{\commI}[1]{\marginpar{%
\begin{color}{blue}
\vskip-\baselineskip 
\raggedright\footnotesize
\itshape\hrule \smallskip I: #1\par\smallskip\hrule\end{color}}}

\newcommand{\commII}[1]{\marginpar{%
\begin{color}{magenta}
\vskip-\baselineskip 
\raggedright\footnotesize
\itshape\hrule \smallskip I: #1\par\smallskip\hrule\end{color}}}




\hyphenation{re-pub-lished}

\parskip 4pt plus 2pt minus 2pt



\def\bfdefault{b}
\overfullrule=5pt

\def \F{{\mathbb F}}
\def \K{{\mathbb K}}
\def \Z{{\mathbb Z}}
\def \Q{{\mathbb Q}}
\def \R{{\mathbb R}}
\def \C{{\\mathbb C}}
\def\Fp{\F_p}
\def \fp{\Fp^*}

\title[Multiple Character Sums]{On Some Multiple Character Sums}

 \author[I. D. Shkredov]{Ilya D. Shkredov}

\address{ Steklov Mathematical Institute of Russian Academy
of Sciences, ul. Gubkina 8, Moscow, Russia, 119991, and Institute for Information Transmission Problems  of Russian Academy
of Sciences, Bolshoy Karet\-ny Per. 19, Moscow, Russia, 127994}
\email{ilya.shkredov@gmail.com}

 \author[I. E. Shparlinski] {Igor E. Shparlinski}

\address{Department of Pure Mathematics, University of New South Wales,
Sydney, NSW 2052, Australia}
\email{igor.shparlinski@unsw.edu.au}

\begin{abstract} We improve a recent result of B. Hanson (2015) on multiplicative character sums
with expressions of the type $a + b +cd$ and variables $a,b,c,d$ from four distinct sets
of a finite field.  We also consider similar sums with $a + b(c+d)$. These bounds  rely on some recent advances in additive combinatorics.
\end{abstract}


\maketitle

\section{Introduction}

\subsection{Motivation and previous results}
Let $p$ be a prime and let $\F_p$ be the finite field of $p$ elements.

%

Given four sets $\cA, \cB, \cC, \cD \subseteq \F_p^*$, and two sequences of  weights
 $\alpha= (\alpha_{a})_{a\in \cA}$, $\beta = \( \beta_{b,c,d}\)_{(b,c,d) \in \cB\times\cC \times \cD}$
 supported on $\cA$ and  $ \cB\times\cC \times \cD$, respectively, we consider the multilinear character sums
\begin{equation}
\label{eq:Sum S}
S_\chi(\cA, \cB, \cC, \cD;\alpha, \beta) =   \sum_{a \in\cA} \sum_{b \in\cB} \sum_{c \in \cC}
 \sum_{d\in \cD}\alpha_{a} \beta_{b,c,d}\chi(a + b + cd) ,
\end{equation}
where $\chi$ is a fixed nontrivial multiplicative character of $\F_p^*$,
see~\cite[Chapter~3]{IwKow} for a background on characters.
We also consider related sums
\begin{equation}
\label{eq:Sum T}
T_\chi(\cA, \cB, \cC, \cD;\alpha, \beta) =   \sum_{a \in\cA} \sum_{b \in\cB} \sum_{c \in \cC}
 \sum_{d\in \cD}\alpha_{a} \beta_{b,c,d}\chi\(a + b(c + d)\).
\end{equation}
In is easy to see that if the weight $\beta_{b,c,d}$ is multilinear, that is of the form
\begin{equation}
\label{eq:MuliLInWeight}
\beta_{b,c,d} = \beta_{b}\gamma_{c}\delta_{d}
 \end{equation}
the sums~\eqref {eq:Sum T}
can easily be reduced to the sums of the form~\eqref{eq:Sum S} (with slightly modified weights).

First we recall that for the bilinear analogues of these sums,
that is, for the sums
$$
S_\chi(\cU, \cV; \varphi, \psi)
=    \sum_{u\in\cU} \sum_{v \in\cV} \varphi_{u} \psi_v \chi(u+v)
$$
with sets $\cU, \cV \subseteq \F_p$ and weights
$\varphi= (\varphi_{u})_{u\in \cU}$, $\psi = \( \psi_{v}\)_{v \in \cV}$,
we have  the classical bound
\begin{equation}
\label{eq:bilin}
|S_\chi(\cU, \cV; \varphi, \psi)| \le  \sqrt{p\varPhi\Psi},
\end{equation}
with
$$
\varPsi =  \sum_{u\in\cU} |\varphi_{u}|^2  \mand
\Psi =  \sum_{v\in\cV} |\psi_{v}|^2.
$$
Indeed, using the trivial inequalities
\begin{align*}
|S_\chi(\cU, \cV; \varphi, \psi)|^2  &\le
  \sum_{u\in\cU} |\varphi_{u}|^2
\sum_{u\in\cU} \left |\sum_{v \in\cV}  \psi_v \chi(u+v)\right |^2\\
& \le   \sum_{u\in\cU} |\varphi_{u}|^2
\sum_{u\in\F_p} \left |\sum_{v \in\cV}  \psi_v \chi(u+v)\right |^2,
\end{align*}
expanding the square, changing the order of summation and recalling
the orthogonality of characters,
we immediately obtain~\eqref{eq:bilin}.

Although our results apply to more general weights, we always assume
that the weighst satisfy the inequalities
\begin{equation}
\label{eq:weight}
\max_{a\in \cA} |\alpha_{a}|\le 1 \mand
\max_{(b,c,d) \in \cB\times\cC \times \cD} | \beta_{b,c,d}|\le 1.
\end{equation}
We also recall that $0$ is excluded from the sets $\cA, \cB, \cC,\cD$
(it is trivial to adjust our bounds to include this case as well).

Thus, using $A$, $B$, $C$ and $D$ to
denote the cardinalities of $\cA, \cB, \cC$ and  $\cD$
respectively, under the condition~\eqref{eq:weight} we easily derive
\begin{equation}
\label{eq:triv}
|S_\chi(\cA, \cB, \cC, \cD; \alpha, \beta)|  \le ABCD \sqrt{\frac{p}{AM}}
\end{equation}
where $M = \max\{B, C, D\}$ and a similar  bound for
$T_\chi(\cA, \cB, \cC, \cD; \alpha, \beta)$ (one only  has to consider the
contribution from the terms with $c+d = 0$ separately).

In the case of constant weights $\alpha_{a} =  \beta_{b,c,d} =1$, these sums,
which we denote as $S_\chi(\cA, \cB, \cC, \cD)$,
have been nontrivially estimated by Hanson~\cite{Han}  under the
 various restrictions on  $A$, $B$, $C$ and $D$.
For example, it is shown in~\cite{Han} that for any fixed $\varepsilon> 0$ there
exists some $\eta >0$ such that if $A, B, C, D  \ge p^\varepsilon$ and
\begin{itemize}
\item either $C < p^{1/2}$ and $A^{56}B^{28}C^{33} D^4 \ge p^{60+\varepsilon}$;
\item or $C  \ge p^{1/2}$ and $A^{112}B^{56} D^8 \ge p^{87+\varepsilon}$;
\end{itemize}
then
\begin{equation}
\label{eq:Hans}
S_\chi(\cA, \cB, \cC, \cD) = O\(ABCD p^{-\eta}\).
\end{equation}
Balog and Wooley~\cite[Section~6]{BalWool} have sugegsted an alternative approach
to bounding the sums $S_\chi(\cA, \cB, \cC, \cD)$ which may lead to more explicit statements
with a small, but explicit values of $\eta$. However our results seems to supersede 
these bounds as well.

Several more bounds of multiplicative character sums, which go beyond
an immediate application of~\eqref{eq:bilin}, have been given by
Bourgain,  Garaev,  Konyagin
and Shparlinski~\cite{BGKS}, Chang~\cite{Chang},
Friedlander and  Iwaniec~\cite{FrIw},
Karatsuba~\cite{Kar}
and Shkredov and Volostnov~\cite{ShkVol}.
However, despite the active interest to such sums of multiplicative characters, there
is a large discrepancy between the strength and generality of the above results and
the results available for similar exponential sums for which strong explicit bounds are
known in very general scenarios,
see~\cite{Bou1, Bou2, BouGar, BouGlib, Gar, Ost, PetShp}.  Here we make a further step towards eliminating this disparity.

\subsection{General notation}

Throughout the paper,  the expressions $A \ll B$,  $B \gg A$ and $A=O(B)$ are each 
equivalent to the statement that $|A|\le cB$ for some positive constant $c$. Throughout the paper,
the implied constants in these symbols may occasionally, where obvious,
 depend on the integer  positive parameter $\nu$, and are absolute otherwise.

\subsection{Main results}

It is convenient to assume that $BCD\le p^2$. Clearly this restriction is not 
important as if it fails we can use $M \ge p^{2/3}$ in the bound~\eqref{eq:triv},
getting a bound of the type~\eqref{eq:Hans} already
for $A \ge p^{1/3+\varepsilon}$ with any fixed $\varepsilon > 0$.

We are now ready to present our main results.

\begin{theorem}
\label{thm:Bound ST}
For any sets  $\cA, \cB, \cC, \cD \subseteq \F_p^*$,
of cardinalities $A,B,C,D$, respectively,
and two sequences of  weights
 $$
 \alpha= (\alpha_{a})_{a\in \cA} \mand
 \beta = \( \beta_{b,c,d}\)_{(b,c,d) \in \cB\times\cC \times \cD},
$$
satisfying~\eqref{eq:weight}, for any fixed integer $\nu \ge 1$,
we have
\begin{equation*}
\begin{split}
S_\chi(\cA, \cB, \cC, \cD&;\alpha, \beta), \ T_\chi(\cA, \cB, \cC, \cD;\alpha, \beta) \\
&  \ll \((BCD)^{1-1/(4\nu)}  +  (BCD)^{1-1/(2\nu)}  M^{1/(2\nu)}\) \\
& \qquad \qquad \qquad \times \left\{
\begin{array}{ll}
A^{1/2}p^{1/2},& \text{if $\nu =1$},\\
 \(A p^{1/(4\nu)} + A^{1/2}p^{1/(2\nu)}\), & \text{if $\nu \ge 2$},
\end{array}
\right.
\end{split}
\end{equation*}
where $M = \max\{B, C, D\}$. 
\end{theorem}

To understand the strength of  Theorem~\ref{thm:Bound ST},
we assume that $A \ge p^{\varepsilon}$ for some $\varepsilon > 0$.
Taking $\nu$ large enough (for example $\nu = \fl{2/\varepsilon} +1$)
and assuming that $M \le (BCD)^{1/2}$ (which certainly holds
in the most interesting case when the sets $\cB$, $\cC$ and
$\cD$ are of comparable sizes), we see that  the bound of  Theorem~\ref{thm:Bound ST}
takes shape
$$
S_\chi(\cA, \cB, \cC, \cD;\alpha, \beta), T_\chi(\cA, \cB, \cC, \cD;\alpha, \beta)   \ll ABCD \(\frac{p}{BCD}\)^{1/(4\nu)},
$$
which is nontrivial as along as $BCD> p^{1+\varepsilon}$ for some $\varepsilon > 0$.

In another interesting case of all sets of asymptotically the same size,
that is, when $A\sim B \sim C \sim D$, taking $\nu =1$ we obtain
$$
S_\chi(\cA, \cB, \cC, \cD;\alpha, \beta),\  T_\chi(\cA, \cB, \cC, \cD;\alpha, \beta) \\
  \ll A^{11/4} p^{1/2}
$$
which is nontrivial for $A\ge  p^{2/5 +\varepsilon}$.

%

\section{Preliminaries}
 \subsection{Background from arithmetic combinatorics}
 \label{sec:AddComb}

For sets $\cB, \cC, \cD \subseteq \F_p^*$.  we denote by  $I(\cB, \cC, \cD)$   the number of solutions to
the equation
$$
b_1 + c_1d_1 = b_2 + c_2d_2 \qquad b_1, b_2 \in \cB, \ c_1,c_2 \in \cC, \ d_1,d_2 \in \cD.
$$
Roche-Newton,   Rudnev and  Shkredov~\cite[Equation~(4)]{RNRS} have shown that the points--planes incidence bound of  Rudnev~\cite{Rud} yields the following estimate:

\begin{lemma}
\label{lem:I Fp}
Let $\cB, \cC, \cD \subseteq \F_p^*$ be of cardinalities $B,C,D$, respectively,
with $BCD \le p^2$.
Then we have
$$
I(\cB, \cC, \cD)  \ll  B^{3/2}C^{3/2}D^{3/2} +  BCD M,
$$
where $M=\max\{B,C,D\}$.
\end{lemma}

Furthermore, for sets $\cB, \cC, \cD \subseteq \F_p^*$.  we denote by  $J(\cB, \cC, \cD)$   the number of solutions to
the equation
$$
b_1(c_1+d_1) = b_2 (c_2+d_2)\qquad  b_1, b_2 \in \cB, \ c_1,c_2 \in \cC, \ d_1,d_2 \in \cD.
$$
It is shown in~\cite[Lemma~2.3]{PetShp} that using
some recent results of
Aksoy Yazici, Murphy, Rudnev and
Shkredov~\cite[Theorem~19]{AYMRS} (which are also based on the work of Rudnev~\cite{Rud}),
one can derive the following analogue of Lemma~\ref{lem:I Fp}

\begin{lemma}
\label{lem:J Fp}
Let $\cB, \cC, \cD \subseteq \F_p^*$  be of cardinalities $B,C,D$, respectively,
with $BCD \le p^2$.  Then we have
$$
J(\cB, \cC, \cD)  \ll  B^{3/2}C^{3/2}D^{3/2} +  BCD M,
$$
where $M=\max\{B,C,D\}$.
\end{lemma}

 \subsection{Bounds of some character  sums on average}
 \label{sec:CharSum}

 The following result is very well-know, and in the case when $\cA$ is an interval
 it dates back to   Davenport \&
Erd{\H o}s~\cite{DavErd}. The proof  transfers to the case of general sets without any changes.
Indeed, for $\nu=1$ it is based on the elementary identity
$$
 \sum_{\lambda \in \F_p}
  \chi\(u+\lambda\)   \overline \chi\(v+\lambda\) =
  \left\{
\begin{array}{ll}
-1,& \text{if $u\ne v$},\\
 p-1, & \text{if $u= v$},
\end{array}
\right.
\qquad u,v  \in \F_p,
$$
where $ \overline \chi$ is the complex conjugate character, see~\cite[Equation~(3.20)]{IwKow}. 
For $\nu \ge 2$, the proof is completely analogous but appeal to 
 the Weil bound of multiplicative character sums, see~\cite[Theorem~11.23]{IwKow}, to estimate 
 ``off-diagonal'' terms.

\begin{lemma}
\label{lem:DavErd} For any set  $\cA  \subseteq \F_p^*$,
of cardinalities $A$
and ta  sequences of  weights
 $ \alpha= (\alpha_{a})_{a\in \cA}$
satisfying~\eqref{eq:weight},   for any fixed integer $\nu \ge 1$,
we have
$$ \sum_{\lambda \in \F_p}
 \left|\sum_{a\in \cA} \alpha_a \chi\(\lambda +a\)
\right|^{2\nu} \ll
 \left\{
\begin{array}{ll}
A p,& \text{if $\nu =1$},\\
 A^{2\nu} p^{1/2} + A^\nu p, & \text{if $\nu \ge 2$}.
\end{array}
\right.
$$
\end{lemma}

\section{Proof of Theorem~\ref{thm:Bound ST}}

 \subsection{Bound on $S_\chi(\cA, \cB, \cC, \cD;\alpha, \beta)$}
 \label{sec:sum S}

 Using~\eqref{eq:weight}, we obtain
$$
|S_\chi(\cA, \cB, \cC, \cD;\alpha, \beta)|\le  \sum_{b \in\cB} \sum_{c \in \cC}
 \sum_{d\in \cD}\left|  \sum_{a \in\cA} \alpha_{a} \chi(a + b + cd)\right| .
 $$
Now, for every $\lambda \in \F_p$ we collect together the terms with the same value of $b  + cd =  \lambda$,
and write
\begin{equation}
\label{eq:S and N}
|S_\chi(\cA, \cB, \cC, \cD;\alpha, \beta)|\le   \sum_{\lambda\in \F_p} K(\cB, \cC, \cD;\lambda)
\left|  \sum_{a \in\cA} \alpha_{a} \chi(a + \lambda )\right|,
\end{equation}
where $K(\cB, \cC, \cD;\lambda)$ is  the number of solutions to
the equation
$$
b  + cd =  \lambda, \qquad (b,c,d) \in \cB \times  \cC \times  \cD.
$$
Clearly we have
\begin{equation}
\begin{split}
\label{eq:I and N}
&\sum_{\lambda \in \F_p} K(\cB, \cC, \cD;\lambda)   = BCD,\\
&\sum_{\lambda \in \F_p} K(\cB, \cC, \cD;\lambda) ^2 = I(\cB, \cC, \cD).
\end{split}
\end{equation}
Therefore, by the H{\"o}lder inequality, for any integer $\nu \ge 1$, we have
\begin{equation*}
\begin{split}
|S_\chi(\cA, \cB&, \cC, \cD;\alpha, \beta)|^{2\nu} \\
& \le
\(\sum_{\lambda \in \F_p} K(\cB, \cC, \cD;\lambda) \)^{2\nu - 2}\\
& \qquad \qquad  \qquad 
\sum_{\lambda \in \F_p} K(\cB, \cC, \cD;\lambda) ^2  \cdot \sum_{\lambda\in \F_p}
\left|  \sum_{a \in\cA} \alpha_{a} \chi(a + \lambda )\right|^{2\nu}\\
& = (BCD)^{2\nu-2} I(\cB, \cC, \cD)\sum_{\lambda\in \F_p}
\left|  \sum_{a \in\cA} \alpha_{a} \chi(a + \lambda )\right|^{2\nu}.
\end{split}
\end{equation*}
Recalling Lemmas~\ref{lem:I Fp} and~\ref{lem:DavErd}, we obtain
\begin{equation*}
\begin{split}
|S_\chi(\cA, \cB, \cC, \cD;\alpha, \beta)|^{2\nu}  &  \ll \((BCD)^{2\nu-1/2}  +  (BCD)^{2\nu-1}  M\)\\
& \qquad \qquad \quad  \times \left\{
\begin{array}{ll}
A p,& \text{if $\nu =1$},\\
 A^{2\nu} p^{1/2} + A^\nu p, & \text{if $\nu \ge 2$},
\end{array}
\right.
\end{split}
\end{equation*}
and the result follows.

\subsection{Bound on $T_\chi(\cA, \cB, \cC, \cD;\alpha, \beta)$}

 We proceed as in the case of the sums $S_\chi(\cA, \cB, \cC, \cD;\alpha, \beta)$.

We define $L(\cB, \cC, \cD;\lambda)$ is  the number of solutions to
the equation
$$
b(c+d) =  \lambda, \qquad (b,c,d) \in \cB \times  \cC \times  \cD, 
$$
and then instead of~\eqref{eq:S and N} write
$$
|T_\chi(\cA, \cB, \cC, \cD;\alpha, \beta)|\le
 \sum_{\lambda\in \F_p} L(\cB, \cC, \cD;\lambda)
\left|  \sum_{a \in\cA} \alpha_{a} \chi(a + \lambda )\right|.
$$
As before, by the  H{\"o}lder inequality, for any integer $\nu \ge 1$, we have
\begin{align*}
|T_\chi(\cA, \cB, \cC, \cD&;\alpha, \beta)|^{2\nu}\\
& \le (BCD)^{2\nu-2} J(\cB, \cC, \cD)\sum_{\lambda\in \F_p}
\left|  \sum_{a \in\cA} \alpha_{a} \chi(a + \lambda )\right|^{2\nu}.
\end{align*}
Using Lemma~\ref{lem:J Fp} instead of Lemma~\ref{lem:I Fp}
in the argument of Section~\eqref{sec:sum S} we obtain
the desired estimate.

\section{Comments}

We note that in the case of the multilinear weights of the
form~\eqref{eq:MuliLInWeight} (an in particular in the case of constant
weights) the roles of the sets $\cA$ and $\cB$ can be interchanged.
Furthermore,  in this case, writing
$$
\chi(a + b + cd)  = \chi(c)  \chi(d + (a + b)c^{-1}) =   \chi(d)  \chi(c + (a + b)d^{-1})
$$
one can obtain the bound of Theorem~\ref{thm:Bound ST}
with any permutation of the roles of $A,B,C,D$.

It is also easy to see that we can abandon the assumption~\eqref{eq:weight}
and obtain a more precise version of  of Theorem~\ref{thm:Bound ST}
with $L^1$ and $L^2$ norms of the weight sequences $\alpha$ and $\beta$.

Finally, using results from \cite{AYMRS} one can obtain similar bounds for several
other character sums exactly in the same way, for example, for the sums
$$
    \sum_{a \in\cA} \sum_{b \in\cB} \sum_{c \in \cC} \sum_{d\in \cD}\alpha_{a} \beta_{b,c,d}  \chi\(a+ \frac{b+c}{d+c}\)\,,  
$$
and 
$$
    \sum_{a \in\cA} \sum_{b \in\cB} \sum_{c \in \cC} \sum_{d\in \cD}\alpha_{a} \beta_{b,c,d} \chi\(a+\frac{b}{c+d}\)\, ,
$$
as well  some other related sums.

We also note that the same approach (with $\nu =1$) applies to the sums
\begin{equation*}
\begin{split}
& U(\cA, \cB, \cC, \cD;\alpha, \beta) =   \sum_{a \in\cA} \sum_{b \in\cB} \sum_{c \in \cC}
 \sum_{d\in \cD}\alpha_{a} \beta_{b,c,d}\ep(a (b+cd)), \\
& V(\cA, \cB, \cC, \cD;\alpha, \beta) =   \sum_{a \in\cA} \sum_{b \in\cB} \sum_{c \in \cC}
 \sum_{d\in \cD}\alpha_{a} \beta_{b,c,d}\ep(ab (c+d)), 
\end{split}
\end{equation*}
where $\ep(u) = \exp(2\pi i u/p)$, 
which are very similar to those appearing in the 
proofs of~\cite[Theorems~1.1 and~1.3]{PetShp}.
Indeed, by the Cauchy inequality, in the notation of Section~\ref{sec:sum S}, we obtain
\begin{equation*}
\begin{split}
|U(\cA, \cB, \cC, \cD;\alpha, \beta)|^{2} 
& \le \sum_{\lambda \in \F_p} K(\cB, \cC, \cD;\lambda) ^2  \cdot \sum_{\lambda\in \F_p}
\left|  \sum_{a \in\cA} \alpha_{a} \ep(a\lambda )\right|^{2}\\
& =   I(\cB, \cC, \cD)
 \sum_{\lambda\in \F_p}
\left|  \sum_{a \in\cA} \alpha_{a} \ep(a\lambda )\right|^{2}.
\end{split}
\end{equation*}
Using the orthogonality of exponential functions, we derive
$$
|U(\cA, \cB, \cC, \cD;\alpha, \beta)|^{2} \le  p  I(\cB, \cC, \cD) A, 
$$
and revoking  Lemma~\ref{lem:I Fp}, we derive 
$$
U(\cA, \cB, \cC, \cD;\alpha, \beta) \ll  p^{1/2} \((B C D)^{3/4} +  (BCD M)^{1/2}\) A^{1/2}, 
$$
where, as before,  $M = \max\{B, C, D\}$. The same argument also gives
$$
V(\cA, \cB, \cC, \cD;\alpha, \beta) \ll  p^{1/2} \((B C D)^{3/4} +  (BCD M)^{1/2}\) A^{1/2}. 
$$
Assuming that $M \le (BCD)^{1/2}$, we see that   
these bounds nontrivial as along as $A(BCD)^{1/2}> p^{1+\varepsilon}$ for some $\varepsilon > 0$.
In particular, under the  condition 
$A\sim B \sim C \sim D$ this holds for $A \ge p^{2/5+\varepsilon}$, which is 
consistent with the range of nontriviality  of~\cite[Theorems~1.1 and 1.3]{PetShp}.

\section*{Acknowledgement}
%
%
During the preparation of this work, 
the second  author
was supported 
by the Australian Research Council Grant DP140100118.

\end{document}